\documentstyle[amssymb,12pt]{amsart}
\newtheorem{apart}{}[section]
\newtheorem{teo}[apart]{{\bf Theorem}}
\newtheorem{prop}[apart]{{\bf Proposition}}
\newtheorem{lem}[apart]{{\bf Lemma}}

\author[de Le\'on ]{\sc M. de Le\'on }
\address[de Le\'on]{Instituto de Ciencias Matem\'aticas (CSIC-UAM-UC3M-UCM), C/ Nicol\'as Ca\-bre\-ra 15,
Campus Cotoblanco UAM, E-28049 Madrid, Spain} \email[de
Le\'on]{mdeleon\symbol{64}icmat.es}

\author[Mart{\'\i}n ]{\sc A. Mart\'{\i}n M\'{e}ndez }
\address[Mart{\'\i}n]{Departamento de Matem\'{a}tica Aplicada II (ETSI
Telecomunicaci\'{o}n), Universidade de Vigo, Lagoas-Marcosende,
36310, Vigo (Pontevedra), Spain}
\email[Mart{\'\i}n]{amartin\symbol{64}dma.uvigo.es}

\title[Principal bundle structures among second order frame bundles]{Principal bundle structures among second order frame bundles}



\date{}

\begin{document}

\maketitle

\begin{abstract}
Using a model for the bundle $\hat{\mathcal F}^2M$ of
semi-holonomic second order frames of a manifold $M$ as an
extension of the bundle ${\mathcal F}^2M$ of holonomic second
order frames of $M$, we introduce in $\hat{\mathcal F}^2M$ a
principal bundle structure over ${\mathcal F}^2M$, the structure
group being the additive group $A_2(n)$ of skew-symmetric bilinear
maps from ${\Bbb R}^2 \times {\Bbb R}^n$ into ${\Bbb R}^n$. The
composition of the projection of that structure with the existing
projection of the bundle $\tilde{\mathcal F}^2M$ of non-holonomic
second order frames of $M$ over $\hat{\mathcal F}^2M$ provides a
principal bundle structure in $\tilde{\mathcal F}^2M$ over
${\mathcal F}^2M$. These results close an existing gap in the
theory of second order frame bundles.
\end{abstract}

\section{Introduction}

\noindent In the recent years, higher order frames have become an
important tool that has contributed to the development of
different chapters both inside the framework of differential
geometry and physics. We can cite for example differential
invariants, classical field theories, quantum field theory,
continuum mechanics~\cite{campos,deleon1,cloud}. Since they were
introduced by Ehresmann in 1995, non-holonomic and semi-holonomic
higher order frame bundles have known several equivalent
interpretations, most of them in terms of jet bundles. The
characterization given here as a starting point is the one
followed in the work by Elzanowski and
Prishepionok~\cite{elzanowski2}, although later we make use of the
interpretation of the non-holonomic and semi-holonomic second
order frame bundles $\tilde{\mathcal F}^2M$ and $\hat{\mathcal
F}^2M$ of a differentiable manifold $M$ as extensions of the
holonomic second order frame bundle ${\mathcal
F}^2M$~\cite{gelines}.

One kind of problems that arises in the context of these theories
is the structure problems. Recently, Brajer\v{c}{\'\i}k, Demko and
Krupka~\cite{krupka} have realized a construction of a principal
bundle structure on the $r$-jet prolongation of the linear frame
bundle of an $n$-dimensional manifold; such prolongations involve
semi-holonomic second order frames~\cite{libermann3}. In this
paper we also deal with a structure problem, raised
in~\cite{paula}, for second order frames: the possible existence
of a principal bundle structure on $\tilde{\mathcal F}^2M$ and
$\hat{\mathcal F}^2M$ over ${\mathcal F}^2M$. The answer to the
question is in the affirmative and the way in which we have
treated it is a constructive way, in the sense that we have
defined in detail such structures. The main drawback lies in
finding a suitable projection $\hat{\pi }^2_2\colon \hat{\mathcal
F}^2M \rightarrow {\mathcal F}^2M$. An early work by
Kol\'a\v{r}~\cite{kolar} provides a symmetrization of
semi-holonomic 2-jets using exact diagrams of vector bundles and
splitting properties. Some other attempt of symmetrization in
higher order using adittional geometric background has only
reached a relative succeed~\cite{mikulski}. The idea underlying in
our definition of the projection $\hat{\pi }^2_2$ lies in the fact
that the description of the bundle $\hat{\mathcal F}^2M$ as an
extension of the bundle ${\mathcal F}^2M$ allows an algebraic
manner to symmetrize the bilinear part of second order jets in a
global way, on the basis that the symmetrization of a
(semi-holonomic) jet somehow commutes with the composition with
holonomic jets. The provided solution also completes the problem
of the existence of principal bundle structures among the three
second order frame bundles $\tilde{\mathcal F}^2M$, $\hat{\mathcal
F}^2M$ and ${\mathcal F}^2M$, the linear frame bundle ${\mathcal
F}M$, and the own manifold $M$. These relations are summarize in
the Section \ref{con} in the final part of the paper.

\section{Preliminars}
\label{pre}

\noindent We are going to deal with five differentiable manifolds:
an $n$-dimensional manifold $M$, its linear frame bundle
${\mathcal F}M$ and the bundles ${\mathcal F}^2M$, $\hat{\mathcal
F}^2M$ and $\tilde{\mathcal F}^2M$ of holonomic, semi-holonomic
and non-holonomic second order frames of $M$, respectively. Let us
denote by $\pi ^1_0$, $\pi ^2_0$, $\hat{\pi }^2_0$ and $\tilde{\pi
}^2_0$ the natural projections of ${\mathcal F}M$, ${\mathcal
F}^2M$, $\hat{\mathcal F}^2M$ and $\tilde{\mathcal F}^2M$ on $M$
and by $G^1(n)=GL(n,{\Bbb R})$, $G^2(n)$, $\hat{G}^2(n)$ and
$\tilde{G}^2(n)$ their respective structure groups. The second
order frames can be described as follows. Let $\tilde{\phi }\colon
{\mathbb R}^n\rightarrow {\mathcal F}M$ be a differentiable map
such that $\pi _0^1 \circ \tilde{\phi }$ is a diffeomorphism. The
1-jet $j_0^1\tilde{\phi }$ with source at the origin $0\in {\Bbb
R}^n$ and target in ${\mathcal F}M$ is a non-holonomic second
order frame at the point $x=\pi _0^1(\tilde{\phi }(0))\in M$. If
$\tilde{\phi }$ satisfies the additional condition $\tilde{\phi
}(0)=j_0^1(\pi _0^1 \circ \tilde{\phi })$, then $j_0^1\tilde{\phi
}$ is called semi-holonomic. A holonomic second order frame is an
invertible $2$-jet $j^2_0f$ where $f$ is a differentiable map from
${\Bbb R}^n$ into $M$ with source at $0\in {\Bbb R}^n$ and target
in $M$; a holonomic second order frame can be also seen as a
semi-holonomic second order frame $j_0^1\tilde{\phi }$ verifying
$\tilde{\phi }=(\pi ^1_0\circ \tilde{\phi } )^{(1)}\circ \eta _1$,
where $\eta _1\colon {\mathbb R}^n \rightarrow {\mathcal
F}{\mathbb R}^n\equiv {\mathbb R}^n \oplus G^1(n)$ is the section
given by $\eta _1={\rm id}_{{\mathbb R}^n}\times I$ and, if $M_1$
and $M_2$ are differentiable manifolds, $F^{(1)}\colon j^1_0g \in
{\mathcal F}M_1 \rightarrow j^1_0(F \circ g)\in {\mathcal F}M_2$
is the prolonged map of a local diffeomorphism $F\colon M_1
\rightarrow M_2$ between the corresponding linear frame bundles.
In order to describe second order frames in local
coordinates~\cite{deleon2,gelines}, we put $\tilde{\phi }
(r^a)\equiv (\phi ^{i}(r^a), \phi _j^{i}(r^a))$, where $\{ r^a \}$
are the natural coordinates in ${\Bbb R}^n$. Then the
non-holonomic frame $j_0^1\tilde{\phi }$ is expressed as
$$\left( \phi ^{i}(0), \phi _j^{i}(0),\frac{\partial \phi
^{i}}{\partial r^j}(0),\frac{\partial \phi _l^{k}}{\partial
r^j}(0)\right)\equiv (x^{i},x^{i}_j,y_j^{i},x_j^{kl}).$$

Here, $(x^{i},x_j^{i})$ are the coordinates of $\tilde {\phi }
(0)$ and $(y_j^{i},x_j^{kl})$ are given by
$$\tilde{\phi }_*(0)\frac{\partial }{\partial
r^j}(0)=y_j^{i}\frac{\partial }{\partial x^{i}}(z)+
x_j^{kl}\frac{\partial }{\partial x_l^{k}}(z),\;\; z=\tilde{\phi
}(0),\;\; j=1,\dots ,n,$$

\noindent with $(y_j^{i})$ non-singular since $\pi _0^1\circ
\tilde{\phi }$ is a diffeomorphism. The semi-holonomic condition
above results in $\phi _j^{i}(0)=\frac{\partial \phi
^{i}}{\partial r^j}(0)$ (i.e., $x_j^{i}=y_j^{i})$ what translates
into a system of local coordinates for $\hat{\mathcal F}^2M$ of
the form $(x^{i},x^{i}_j,x_j^{kl})$. On the other hand, the
holonomic condition leads to $\frac{\partial \phi _l^{k}}{\partial
r^j}(0)=\frac{\partial ^2\phi ^k}{\partial r^l\partial r^j}(0)$,
i.e., $x_j^{kl}=x_l^{kj}$.

Information of another models for non-holonomic and semi-holonomic
second order frames can be found, for example,
in~\cite{deleon1,alberto,paula,gelines,yuen}.

Let $L_2(n)$ be the additive group of bilinear maps from ${\mathbb
R}^n\times {\mathbb R}^n$ into ${\mathbb R}^n$ and let $S_2(n)$
and $A_2(n) $ be the additive subgroups of symmetric bilinear maps
and skew-symmetric bilinear maps from ${\mathbb R}^n\times
{\mathbb R}^n$ into ${\mathbb R}^n$. Let us denote by $\{ E_{i}
\}$ the canonical basis in ${\Bbb R}^n$. If $f\in L_2(n)$ we
define the bilinear maps $f^t$, $f_s$, $f_a \colon {\Bbb R}^n
\times {\Bbb R}^n \rightarrow {\Bbb R}^n$ by
$f^t(E_{i},E_j)=f(E_j,E_{i})$, $f_s=\frac{f+f^t}{2}$ and
$f_a=\frac{f-f^t}{2}$; then $f_s\in S_2(n)$, $f_a\in A_2(n)$ and
$f=f_s+f_a$. The groups $\tilde{G}^2(n)$, $\hat{G}^2(n)$ and
$G^2(n)$ can be represented as the semi-direct product
\begin{eqnarray*}
\tilde{G}^2(n) & \equiv & G^1(n)\times G^1(n)\times L_2(n), \\
\hat{G}^2(n) & \equiv & G^1(n)\times L_2(n), \\
G^2(n) & \equiv & G^1(n)\times S_2(n),
\end{eqnarray*}

\noindent with the group operations given by
$$(a,b,f)(a^{\prime
},b^{\prime },f^{\prime })  =  (aa^{\prime }, bb^{\prime }, a\circ
f^{\prime }+f(a^{\prime },b^{\prime }))$$

\noindent for $\tilde{G}^2(n)$ and
$$(a,f)(a^{\prime },f^{\prime
})  =  (aa^{\prime }, a\circ f^{\prime } + f(a^{\prime },a^{\prime
}))$$

\noindent for $\hat{G}^2(n)$ and $G^2(n)$. For later use we state
the following result, whose proof is straightforward.
\begin{lem}
\label{prel1} If $a\in G^1(n)$ and $f\in L_2(n)$, then $(a\circ
f)^t=a\circ f^t$, $(a\circ f)_s=a\circ f_s$, $(a\circ f)_a=a\circ
f_a$ and $[f(a,a)]^t=f^t(a,a)$, $(f(a,a))_s=f_s(a,a)$,
$(f(a,a))_a=f_a(a,a)$.
\end{lem}

In the last section we make an abridgement of the principal bundle
relations among the five manifolds mentioned at the beginning of
this section. Consequently, it is pertinent to recall that
$\tilde{\mathcal F}^2M$, $\hat{\mathcal F}^2M$ and ${\mathcal
F}^2M$ are all principal bundles over ${\mathcal F}M$ with
projections $\tilde{\pi }^2_1$, $\hat{\pi }^2_1$ and $\pi ^2_1$
carrying a jet $j_0^1\tilde{\phi }$, with $\tilde{\phi }\colon
{\mathbb R}^n\rightarrow {\mathcal F}M,$ into $\tilde{\phi }(0)$.
The respective structure groups are $\tilde{G}^2_1(n)$, $L_2(n)$
and $S_2(n)$. The group $\tilde{G}^2_1(n)$ is
$\tilde{G}^2_1(n)=G^1(n)\times L_2(n)$ endowed with the law
$(a,f)(a^{\prime },f^{\prime })=(aa^{\prime },f^{\prime} +f (I,
a^{\prime }))$~\cite{paula}.

\section{The groups $\hat{G}^2(n)$, $G^2(n)$ and $A_2(n)$}
\label{gro}

\noindent In this section we are going to study some relations
between the groups $\hat{G}^2(n)$, $G^2(n)$ and $A_2(n)$. First of
all, let us note that the inverse $(a,f)^{-1}$ of $(a,f)\in
\hat{G}^2(n)$ is $(a,f)^{-1}=(a^{-1},-a^{-1}\circ
f(a^{-1},a^{-1}))$; if $(b,g)\in \hat{G}^2(n)$, the conjugation of
$(b,g)$ by $(a,f)$ is given by

\begin{eqnarray*}
&  & (a,f)\cdot (c,g)\cdot (a,f)^{-1}=(a,f)\cdot (b,g)\cdot
(a^{-1},-a^{-1}\circ f(a^{-1},a^{-1})) \\
 & = & (ab,a \circ g + f(b,b))\cdot (a^{-1},-a^{-1}\circ
f(a^{-1},a^{-1})) \\
 & = & (aba^{-1},ab \circ [-a^{-1} \circ f(a^{-1},a^{-1})]+((a \circ
g)+f(b,b))(a^{-1},a^{-1}))) \\
 & = & (aba^{-1},-a \circ b \circ a^{-1} \circ f (a^{-1},a^{-1}) + a
\circ g (a^{-1},a^{-1}) + f(ba^{-1},ba^{-1})).
\end{eqnarray*}

In particular, if $b=I$ we obtain
$$(a,f)\cdot (I,g)\cdot (a,f)^{-1}=(I,-f(a^{-1},a^{-1})+a \circ g
(a^{-1},a^{-1}) + f(a^{-1},a^{-1}))=(I,a \circ g(a^{-1},a^{-1}))$$

\noindent and the conjugation of the element $(I,g)\in
\hat{G}^2(n)$ by $(a,f)\in \hat{G}^2(n)$ does not depend on $f$
(~\cite{deleon3}).

\begin{lem}
\label{grol1} Let $(a,f)\in \hat{G}^2(n)$. \begin{itemize}
\item[i)] If $(I,h)\in G^2(n)$, the conjugation $(a,f)\cdot (I,h)
\cdot (a,f)^{-1}$ belongs to $G^2(n)$. \item[ii)] Let us consider
the subgroup of $\hat{G}^2(n)$ defined by the semi-direct product
$G^1(n)\times A_2(n)$. If $(I,h)\in G^1(n)\times A_2(n)$, then the
conjugation $(a,f)\cdot (I,h) \cdot (a,f)^{-1}$ belongs to
$G^1(n)\times A_2(n)$. \item[iii)] There exist unique $(b,g)\in
G^2(n)$ and $(I,h)\in A_2(n)$ such that $(a,f)=(b,g)(I,h)$.
\end{itemize}
\end{lem}

{\bf Proof:}

$ $

Let us put
$$a(E_{i})=a^j_{i}E_j\;\;\; a^{-1}(E_{i})=b^j_{i}E_j \;\;\;
h(E_{i},E_j)=h^{li}_jE_l.$$

It follows that
\begin{eqnarray*}
(a \circ h(a^{-1},a^{-1}))(E_{i},E_j) & = & a \circ
h(a^{-1}E_{i},a^{-1}E_j) \\
=a \circ h (b^r_{i}E_r,b^s_jE_s) & = & a \circ (b^r_{i}b^s_j
h(E_r,E_s)) \\
=a \big( b^r_{i}b^s_jh^{kr}_sE_k \big) & = & b^r_{i}b^s_jh^{kr}_s
a(E_k)=b^r_{i}b^s_jh^{kr}_sa^l_kE_l \end {eqnarray*}

\noindent and if $a \circ h
(a^{-1},a^{-1})(E_{i},E_j)=m^{li}_jE_l$, then
$$m^{li}_j=b^r_{i}b^s_jh^{kr}_sa^l_k.$$

Let us suppose that $(I,h)\in G^2(n)$, i.e., $h^{kr}_s=h^{ks}_r$.
Then
$$m^{lj}_{i}=b^r_jb^s_{i}h^{kr}_sa^l_k=b^s_{i}b^r_jh^{kr}_sa^l_k=b^s_{i}b^r_jh^{ks}_ra^l_k=b^r_{i}b^s_jh^{kr}_sa^l_k=m^{li}_j$$

\noindent (note that in the second last equality we have
interchanged the indexes $r$ and $s$). That means that $a\circ
h(a^{-1},a^{-1})\in G^2(n)$. This proves i). Statement ii) can be
proved in the same way using that $(a,h)\in G^1(n)\times A_2(n) $
if and only if $h^{ij}_k = - h^{ik}_j.$ In order to prove iii),
only note that the elements $(b,g)$ e $(I,h)$ are $(b,g)=(a,f_s)$
and $(I,h)=(I,a^{-1}\circ f_a)$. \\

Statement iii) in Lemma \ref{grol1} points out that $\hat{G}^2(n)$
can be also read as the semi-direct product $\hat{G}^2(n) \equiv
G^2(n) \times A_2(n)$ in such a way that the law group in $G^2(n)
\times A_2(n)$ springs from jet composition.

From now on, we shall often identify $(I,h)$ with $h$; with this
identification we obtain that $S_2(n)$ and $A_2(n)$ are closed
subgroups of $\hat{G}^2(n)$~\cite{deleon3} and

\begin{lem}
\label{grol3} The additive groups $S_2(n)$ and $A_2(n)$ are normal
subgroups of $\hat{G}^2(n)$.
\end{lem}

Since $A_2(n)$ is a normal closed subgroup of $\hat{G}^2(n)$, then
the quotient group $\frac{\hat{G}^2(n)}{A_2(n)}$ is a Lie group.

\begin{prop}
\label{grop1} $\frac{\hat{G}^2(n)}{A_2(n)}$ and $G^2(n)$ are
isomorphic Lie groups.
\end{prop}

{\bf Proof:}

$ $

Define
$$\mu  \colon  (a,f)A_2(n) \in \frac{\hat{G}^2(n)}{A_2(n)}  \rightarrow  (a,f_s) \in G^2(n). $$

If $(a,f)A_2(n)=(a^{\prime },f^{\prime })A_2(n)$, there exists
$(I,h)\in A_2(n)$ such that
\begin{eqnarray*} & & (a,f)(a^{\prime },f^{\prime })^{-1}=(a,f)(a^{\prime
-1},-a^{\prime -1}\circ f^{\prime } (a^{\prime -1},a^{\prime -1}))
\\
 & = & (aa^{\prime -1},-a\circ a^{\prime -1}\circ f^{\prime }
(a^{\prime -1},a^{\prime -1}) + f(a^{\prime -1},a^{\prime -1})).
\end{eqnarray*}

It follows then that $a=a^{\prime }$ and $-f^{\prime }
(a^{-1},a^{-1}) + f(a^{-1},a^{-1})\in A_2(n)$. Since $a^{-1}$
represents an isomorphism, the last expression means that
$-f^{\prime }+f$ belongs to $ A_2(n)$. Therefore $(-f^{\prime
}+f)^t=-(-f^{\prime }+f)=f^{\prime }-f$; hence $ f^t+f=f^{\prime
t}+f^{\prime } $ and $ f_s=f^{\prime }_s$. So we have obtained
$\mu ((a,f)A_2(n))=\mu ((a^{\prime },f^{\prime })A_2(n))$ and $\mu
$ is a well defined map. Moreover,
\begin{itemize} \item[$\bullet $ ] $(a,f)A_2(n)\cdot (a^{\prime
},f^{\prime })A_2(n)=(a,f) (a^{\prime },f^{\prime
})A_2(n)=(aa^{\prime }, a\circ f^{\prime }+ f(a^{\prime },a^{\prime
}))A_2(n);$

\item[$\bullet $ ] $\mu ((a,f)A_2(n)\cdot (a^{\prime },f^{\prime
})A_2(n))=(aa^{\prime }, (a\circ f^{\prime }+ f(a^{\prime }
,a^{\prime })_s);$

\item[$\bullet $ ] $\mu ((a,f)A_2(n)) \cdot \mu ((a^{\prime
},f^{\prime })A_2(n))=(a,f_s)(a^{\prime },f^{\prime }_s)=(aa^{\prime
},a\circ f^{\prime }_s+ f_s(a^{\prime },a^{\prime })).$
\end{itemize}

Using Lemma \ref{prel1} we obtain $a\circ f^{\prime }_s+
f_s(a^{\prime },a^{\prime })=(a\circ f^{\prime }+ f(a^{\prime }
,a^{\prime })_s$ and $\mu $ is a group homomorphism. In order to
prove that $\mu $ is injective, put $\mu ((a,f)A_2(n))=\mu
((a^{\prime },f^{\prime })A_2(n))$. Then, $a=a^{\prime } \; {\rm
y}$ and $f_s=f^{\prime }_s$. According to the calculation
developed above to show that $\mu $ is a well-defined map, it
follows that proving $(a,f)A_2(n)=(a^{\prime },f^{\prime })A_2(n)$
is the same as proving $aa^{\prime -1} = I$ (what is true since
$a=a^{\prime }$) and $-f^{\prime } (a^{-1},a^{-1}) +
f(a^{-1},a^{-1})\in A_2(n)$, i.e., $(-f^{\prime
}+f)(a^{-1},a^{-1})\in A_2(n)$, or,
$$(-f^{\prime
}+f)(a^{-1},a^{-1})(E_{i},E_j)=(f^{\prime
}-f)(a^{-1},a^{-1})(E_j,E_{i}).$$

Now then, since $f_s=f^{\prime }_s$ we obtain $-f^{\prime
}+f=-f^{\prime }_s -f^{\prime }_a+f_s+f_a=-f^{\prime }_a+f_a$, and
therefore,
$$\begin{array}{ccccc}
 &  & (-f^{\prime }+f)(a^{-1},a^{-1})(E_{i},E_j) & = & (-f^{\prime
}_a+f_a)(a^{-1},a^{-1})(E_{i},E_j)\\
 & = & (-f^{\prime }_a+f_a)(a^{-1}(E_{i}),a^{-1}(E_j)) & = & (f^{\prime
}_a-f_a)(a^{-1}(E_{j}),a^{-1}(E_{i}))\\
 & = & (f^{\prime }_a-f_a)(a^{-1},a^{-1})(E_{j},E_{i}) & = & (f^{\prime
}-f)(a^{-1},a^{-1})(E_{j},E_{i}),
\end{array}$$

\noindent and $\mu $ is injective. Finally, given $(b,g)\in
G^2(n)$ and since $g=g_s$, there exists $(b,g)A_2(n)\in
\frac{\hat{G}^2(n)}{A_2(n)}$ satisfying $\mu
((b,g)A_2(n))=(b,g_s)=(b,g)$, and $\mu $ is surjective.\\

There are several group structures on $G^1(n)\times L_2(n)$ in the
literature. For example, the following two laws

\begin{eqnarray*}
(a,f)(a^{\prime },f^{\prime }) & = & (aa^{\prime },a^{\prime -1
}\circ f(a^{\prime },a^{\prime
})+f^{\prime })\\
(a,f)(a^{\prime },f^{\prime }) & = & (aa^{\prime },f+a\circ
f^{\prime }(a^{-1},a^{-1})),
\end{eqnarray*}

\noindent give both of them rise to groups which are isomorphic
with $\hat{G}^2(n)$~\cite{deleon2}.

The paper~\cite{krupka} by Brajer\v{c}{\'\i}k, Demko and Krupka
deals with principal bundle structures in jet prolongations of
higher order which include second order frames. There, the base
manifold of the bundles is the manifold $M$.
Libermann~\cite{libermann3} proved that the first order jet
prolongation $J^1{\mathcal F}M$ of the linear frame bundle
${\mathcal F}M$, which consists of $1$-jets of local sections of
${\mathcal F}M$, can be identified with $\hat{\mathcal F}^2M$. The
Lie group named $(T^r_nL^1_n,\ast )$ in~\cite{krupka} is the
structure group of the principal bundle structure given there for
the $r$-order jet prolongation $J^r{\mathcal F}M$. So it is not
surprising that for $r=1$ such structure group is isomorphic with
$\hat{G}^2(n)$. The group law in $(T^1_nL^1_n,\ast )$ is described
in~\cite{krupka} in local coordinates as
$$(a^{i}_j,a^{ij}_k)(c^{i}_j,c^{ij}_k)=(a^{i}_mc^m_j,a^{il}_kc^l_j+a^{i}_lc^{lj}_mb^m_k),$$

\noindent where $(b^{i}_j)=(a^{i}_j)^{-1}$, and we have

\begin{lem}
\label{grol4} The group law $\ast $ in $T^1_nL^1_n$ can be
expressed as
$$(a,f)(a^{\prime },f^{\prime
})=(aa^{\prime },f (a^{\prime }, I)+a\circ f^{\prime } (I,
a^{-1}))$$ and the group $(T^1_nL^1_n,\ast )$ is isomorphic with
$\hat{G}^2(n)$ through the isomorphism of Lie groups
$$\tau \colon (a,f)\in T^1_nL^1_n \rightarrow \tau (a,f)=(a,f(I,a))\in \hat{G}^2(n).$$
\end{lem}

{\bf Proof:}

$ $

If $a(E_{i})=a^j_{i}E_j$, $a^{-1}(E_{i})=b^j_{i}E_j $, $
f(E_{i},E_j)=a^{li}_jE_l$, $a^{\prime }(E_{i})=a^j_{i}E_j$ and
$f^{\prime }(E_{i},E_j)=c^{li}_jE_l$, the description in local
coordinates for $f (a^{\prime }, I)+a\circ f^{\prime } (I,
a^{-1})$ is given by

\begin{eqnarray*}
 & & [f (a^{\prime }, I)+a\circ f^{\prime } (I,
a^{-1})](E_j,E_k) = f (a^{\prime }, I)(E_j,E_k)+a\circ
f^{\prime } (I, a^{-1})(E_j,E_k)\\
 & = & f(a^{\prime }(E_j),E_k)+a\circ f^{\prime
}(E_j,a^{-1}(E_k))  = f(c^{i}_jE_{i},E_k)+a\circ f^{\prime
}(E_j,b^{i}_kE_{i})\\
 & = & c^{i}_jf(E_{i},E_k)+a\circ (b^{i}_kf^{\prime
}(E_j,E_{i}))  =  c^{i}_ja^{ri}_kE_r+a\circ
(b^{i}_kc^{lj}_{i}E_l)\\
 & = & c^{i}_ja^{ri}_kE_r+b^{i}_kc^{lj}_{i}a(E_l)  =
 c^{i}_ja^{ri}_kE_r+b^{i}_kc^{lj}_{i}a^r_lE_r\\
 & = & (c^r_ja^{ir}_k+b^r_kc^{lj}_ra^{i}_l)E_{i}
\end{eqnarray*}

That means that the group law $\ast $ works as stated. Moreover we
have

\begin{eqnarray*}
& \bullet  & \tau ((a,f)(a^{\prime },f^{\prime }))  =  \tau
((aa^{\prime },f (a^{\prime }, I)+a\circ f^{\prime }
(I, a^{-1})))\\
 & = & (aa^{\prime },f (a^{\prime }, I)(I,aa^{\prime
})+a\circ f^{\prime } (I a^{-1})(I,aa^{\prime }))  = (aa^{\prime
},f(a^{\prime },aa^{\prime })+a\circ f^{\prime }(I,a^{\prime }))\\
& \bullet  & \tau (a,f)\tau (a^{\prime },f^{\prime
})=(a,f(I,a))(a^{\prime },f^{\prime }(I,a^{\prime }))  =
(aa^{\prime },a\circ f^{\prime }(I,a^{\prime })+f(I,a)(a^{\prime
},a^{\prime })) \\
 & = &(aa^{\prime },a\circ f^{\prime }(I,a^{\prime })+f(a^{\prime
},aa^{\prime })) \end{eqnarray*}

\noindent and $\tau $ is an isomorphism. Of course, $\tau
^{-1}(a,f)=(a,f(I,a^{-1}))$, and $\tau $ and $\tau ^{-1}$ allow to
recover the group law in $T^1_nL^1_n$ making $(a,f)(a^{\prime
},f^{\prime })=\tau ^{-1}(\tau (a,f) \tau (a^{\prime },f^{\prime
}))$

\section{Remaining bundle structures}
\label{rbs}

\noindent Let ${\mathcal F}^2M^{\hat{G}^2(n)}={\mathcal F}^2M
\times _{G^2(n)}\hat{G}^2(n)=\frac{ {\mathcal F}^2M \times
\hat{G}^2(n)}{ \sim }$ be the quotient manifold obtained from
${\mathcal F}^2M \times \hat{G}^2(n)$ and the equivalence relation
$$(p,k) \sim (p^{\prime },k^{\prime }) \; \Leftrightarrow \; {\rm
there} \; {\rm exists} \; \alpha \in G^2(n) \; {\rm such} \; {\rm
that} \; p^{\prime }=p\alpha ,\; k^{\prime }= \alpha ^{-1}k.$$

The bundle $[(p,k)] \in {\mathcal F}^2M^{\hat{G}^2(n)} \rightarrow
\pi ^2_0(p) \in M$ is called the extension of ${\mathcal F}^2M$ by
the group $\hat{G}^2(n)$. The map
$$\vartheta \colon [(p,k)] \in {\mathcal F}^2M^{\hat{G}^2(n)}
\rightarrow \vartheta ([(p,k)])=pk \in \hat{\mathcal F}^2M$$

\noindent is an isomorphism of principal bundles ~\cite{gelines}.

Let us define a projection $\hat{\pi }^2_2 \colon \hat{\mathcal
F}^2M \rightarrow {\mathcal F}^2M$ as
$$\hat{\pi }^2_2(pk)=\hat{\pi }^2_2 (\vartheta ([(p,k)]))=p \cdot
(a,f_s), \;\; {\rm for} \; k=(a,f).$$

Before seeing that $\hat{\pi }^2_2$ is a well-defined map, we pose
the following basic result, which shows the way in what
composition with semi-holonomic jets and symmetrization commutes.

\begin{prop}
\label{rbsp1} Let $(a,f)\in \hat{G}^2(n)$ and $(b,g)\in G^2(n)$.
If $(b,g)(a,f)=(c,h)$, then $(b,g)(a,f_s)=(c,h_s)$.
\end{prop}

{\bf Proof:}

$ $

Bearing in mind that $(b,g)(a,f)=(ba,b\circ f + g(a,a))$ and
$(b,g)(a,f_s)=(ba,b\circ f_s + g(a,a))$, we have $h=b\circ f +
g(a,a)$. The fact that $h_s=b\circ f_s + g(a,a)$ follows from
$g=g^t$ and Lemma \ref{prel1}.\\

In order to see that $\hat{\pi }^2_2$ is a well-defined map, let
$(p^{\prime },k^{\prime }) \sim (p,k)$. Let $\alpha =(b,g)\in
G^2(n)$ such that $p^{\prime }=p\alpha ,\; k^{\prime }= \alpha
^{-1}k$. We have that
\begin{eqnarray*}
& & k^{\prime }=\alpha ^{-1}k=(b,g)^{-1}(a,f)=(b^{-1},-b^{-1}\circ
g(b^{-1},b^{-1}))(a,f)\\
& =& (b^{-1}a,b^{-1}\circ f+(-b^{-1}\circ
g(b^{-1},b^{-1}))(a,a))=(b^{-1}a,b^{-1}\circ
(f-g(b^{-1}a,b^{-1}a))).
\end{eqnarray*}

Let us put $h=b^{-1}\circ (f-g(b^{-1}a,b^{-1}a))$. We must prove
that
$$\hat{\pi }^2_2 ([(p^{\prime },k^{\prime })])=\hat{\pi }^2_2 ([(p\alpha ,\alpha
^{-1}k)])=p\alpha \cdot (b^{-1}a,h_s)=p \cdot (a,f_s)=\hat{\pi
}^2_2 ([(p,k)]).$$

Since
$$(p\alpha )\cdot (b^{-1}a,h_s)=(p \cdot (b,g))\cdot
(b^{-1}a,h_s)=p \cdot ((b,g)\cdot (b^{-1}a,h_s))$$

\noindent it is enough to see that $(b,g)\cdot
(b^{-1}a,h_s)=(a,f_s)$; but this follows from Proposition
\ref{rbsp1}. So, $\hat{\pi }^2_2 $ is a well-defined map.
Moreover, $\hat{\pi }^2_2 $ is surjective since given $q\in
{\mathcal F}^2M$ we have $\hat{\pi }^2_2(q)=\hat{\pi
}^2_2(q(I,0))=q\cdot (I,0)=q$.

\begin{lem}
\label{rbsl1} $\hat{\pi }^2_0=\pi ^2_0 \circ \hat{\pi }^2_2$.
\end{lem}

{\bf Proof:}

$ $

Let us consider $pk \in \hat{\mathcal F}^2M$, with $p\in {\mathcal
F}^2M$ and $k=(a,f)\in \hat{G}^2(n)$. Then
$$\hat{\pi }^2_0 (pk)=\hat{\pi }^2_0 (p)=\pi ^2_0 (p)=\pi
^2_0(p(a,f_s))=\pi ^2_0(\hat{\pi }^2_2 (pk)),$$ \noindent since
$k\in \hat{G}^2(n)$, $(a,f_s) \in G^2(n)$ and the fibers over
$\hat{\pi }^2_0 (p)=\pi ^2_0 (p)$ in the bundles $\hat{\mathcal
F}^2M \rightarrow M$ and ${\mathcal F}^2M \rightarrow M$ are,
respectively, $p\cdot \hat{G}^2(n)$ and $p\cdot G^2(n)$.

\begin{lem}
\label{rbsl2} Given $q\in {\mathcal F}^2M$ we have $(\hat{\pi
}^2_2)^{-1}(q)=q \cdot A_2(n)$.
\end{lem}

{\bf Proof:}

$ $

" $\subset$ "  Let $p\in {\mathcal F}^2M$ and $k\in \hat{G}^2(n)$
be such that $pk \in (\hat{\pi }^2_2)^{-1}(q)$. Using Lemma
\ref{rbsl1} we obtain $\hat{\pi }^2_0 (pk)=\pi ^2_0(\hat{\pi }^2_2
(pk))=\pi ^2_0(q)$. Therefore, $pk$ and $q$ lie in the same fiber
in $\hat{\pi }^2_0\colon \hat{\mathcal F}^2M \rightarrow M$ and
there exists $\tilde{k}=(\tilde{a},\tilde{f})\in \hat{G}^2(n)$
verifying $pk=q\tilde{k}$. Consequently,
$$q(I,0)=q=\hat{\pi }^2_2 (pk)=\hat{\pi }^2_2 (q\tilde{k})=q\cdot
(\tilde{a},\tilde{f}_s)).$$

Since the action is free, we have $\tilde{a}=I$ and
$\tilde{f}_s=0$, i.e., $\tilde{a}=I$ and $\tilde{f}=-\tilde{f}^t$.
Therefore, $(I,\tilde{f})\in A_2(n)$ and we have proved that $pk$
can be written in the form $q(I,\tilde{f})$ with
$\tilde{f}=-\tilde{f}^t$ that is $pk\in q\cdot A_2(n)$.

" $\supset $ "  Let us consider $h\in A_2(n)$. From the definition
of $\hat{\pi }^2_2$ we obtain, since $h^t=-h$, that $\hat{\pi
}^2_2 (q(I,h))=q(I,\frac{h+h^t}{2})=q(I,0)=q$.

\begin{lem}
\label{rbsl3} $\hat{\pi }^2_1=\pi ^2_1 \circ \hat{\pi }^2_2$.
\end{lem}

{\bf Proof:}

$ $

Let us consider $pk\in \hat{\mathcal F}^2M$, with $p\in {\mathcal
F}^2M$ and $k\in \hat{G}^2(n)$. Let us put $\hat{\pi }^2_2
(pk)=q$. Therefore we can write $pk=q(I,h)$ with $(I,h)\equiv h\in
A_2(n)$. Since $\hat{\pi }^2_1\colon \hat{\mathcal F}^2M
\rightarrow {\mathcal F}M$ is a principal bundle with $L_2(n)$ as
structure group and $A_2(n)\subset L_2(n)$, then
$$\hat{\pi }^2_1 (q(I,h))=\hat{\pi }^2_1 (q)=\pi ^2_1(q),$$

\noindent since $\pi ^2_1 = (\hat{\pi }^2_1)_{|{\mathcal F}^2M}$.
Hence, $\hat{\pi }^2_1(pk)=\hat{\pi }^2_1 (q(I,h))=\pi ^2_1(q)=\pi
^2_1 ( \hat{\pi }^2_2 (pk))$.\\

Let us also note that the projection $\hat{\pi }^2_2$ is
differentiable since so is the map $f\in L_2(n) \rightarrow f_s
\in S_2(n)$.

\begin{lem}
\label{rbsl5} Let $G$ be a Lie group. Let $H_1$ and $H_2$ be Lie
subgroups of $G$ verifying that for every element $g\in G$ there
exist unique $h_1\in H_1$, $h_2\in H_2$ such that $g=h_1h_2$. Let
$M$ and $N$ be differentiable manifold. Let $\hat{\xi } \colon
M\rightarrow G$, $\xi \colon N\rightarrow H_1$ and $f\colon
M\rightarrow N$ be differentiable maps. The map $\sigma \colon M
\rightarrow H_2$ defined by $\hat{\xi }(x)=\xi (f(x)) \sigma (x)$
is differentiable.
\end{lem}

{\bf Proof:}

$ $

It is enough to notice that $\sigma (x)=(\xi (f(x)))^{-1} \hat{\xi
}(x)$.

\begin{teo}
\label{rbst1} $\hat{\mathcal F}^2M$ is a principal bundle over
${\mathcal F}^2M$ with structure group $A_2(n)$ and projection
$\hat{\pi }^2_2$.
\end{teo}

{\bf Proof:}

$ $

The action on the right of the group $A_2(n)$ on $\hat{\mathcal
F}^2M$ the we are going to consider is the restriction of the
action on the right of $\hat{G}^2(n)$ on $\hat{\mathcal F}^2M$;
so, such action is free.

We shall identify $\hat{\mathcal F}^2M$ with the image of
${\mathcal F}^2M^{\hat{G}^2(n)}$ by the map $\vartheta $ given
above. Let us define
$$\Omega  \colon pkA_2(n) \in \frac{\hat{\mathcal F}^2M}{A_2(n)}
\rightarrow \Omega (pkA_2(n))=p(a,f_s)=\hat{\pi }^2_2 (pk)
{\mathcal F}^2M $$

\noindent where $k=(a,f)\in \hat{G}^2(n)$.

The definition of $\Omega $ does not depend on the element on the
equivalence class $pkA_2(n)=p(a,f)A_2(n)=\{(p(a,f))(I,g) \;;\;
g\in A_2(n)\}$. Indeed, $(p(a,f))(I,g)=p((a,f)(I,g))=p(a, a\circ g
+f)$; so, using Lemma \ref{prel1} and $g=-g^t$, we have
\begin{eqnarray*}
& \Omega  & (pk(I,g)A_2(n))=\Omega (p(a, a\circ g +f)A_2(n)) \\
 & = & p\left( [a\circ g +f]_s\right)
=p(a,f_s)=\Omega (pkA_2(n))
\end{eqnarray*}

Now, let $pkA_2(n)$, $p^{\prime }k^{\prime }A_2(n) \in
\frac{\hat{\mathcal F}^2M}{A_2(n)}$ be such that $\Omega (pk
A_2(n))=\Omega (p^{\prime }k^{\prime }A_2(n)$. Then $\hat{\pi
}^2_2(pk)=\hat{\pi }^2_2(p^{\prime }k^{\prime })$, so that $pk$
and $p^{\prime }k^{\prime }$ lie in the same fiber of the
projection $\hat{\pi }^2_2 \colon \hat{\mathcal F}^2M \rightarrow
{\mathcal F}^2M$; by Lemma \ref{rbsl2} there exists $g\in A_2(n)$
such that $pk(I,g)=p^{\prime }k^{\prime }$, i.e.,
$pkA_2(n)=p^{\prime }k^{\prime }A_2(n)$ and $\Omega $ is
injective. $\Omega $ is also surjective since, given $q\in
{\mathcal F}^2M$ there exists $qA_2(n) \in \frac{\hat{\mathcal
F}^2M}{A_2(n)} $ such that $\Omega (qA_2(n))=\hat{\pi }^2_2(q)=q$.
Let us finally note that $\Omega $ and $\Omega ^{-1}$ are
differentiable since so is $\hat{\pi }^2_2$.

Hence we have proved that ${\mathcal F}^2M$ is the quotient space
of $\hat{\mathcal F}^2M$ by the equivalence relation induced by
$A_2(n)$. Moreover, if $\pi \colon p\rightarrow \pi(p)=[p]$ is the
canonical projection, then $\Omega \circ \pi =\hat{\pi }^2_2$.

Next, let $V\subset M$ be an open set and let $\Xi \colon (\pi
^2_0)^{-1}(V) \rightarrow V \times G^2(n)$ and $\hat{\Xi } \colon
(\hat{\pi }^2_0)^{-1}(V) \rightarrow V \times \hat{G}^2(n)$ be
trivializations for the bundles $\pi ^2_0 \colon {\mathcal F}^2M
\rightarrow M$ and $\hat{\pi } ^2_0 \colon \hat{\mathcal F}^2M
\rightarrow M$, respectively, that we choose in such a way that
$\hat{\Xi }(q)=\Xi (q)$ for every $q\in (\pi ^2_0)^{-1}(V)$. Let
us put $U=(\pi ^2_0)^{-1}(V)$. From Lemma \ref{rbsl1} we obtain
$(\hat{\pi }^2_0)^{-1}(V)=(\hat{\pi }^2_2)^{-1}(\pi
^2_0)^{-1}(V)=(\hat{\pi }^2_2)^{-1}(U)$. Define
$$\Sigma \colon p\in (\hat{\pi }^2_2)^{-1}(U) \rightarrow \Sigma (p)
=(\hat{\pi }^2_2(p),\sigma (p)) \in U\times A_2(n) $$

\noindent where $\sigma \colon (\hat{\pi }^2_2)^{-1}(U)
\rightarrow A_2(n)$ is given as follows: if $\Xi (q)=(\pi
^2_0(q),\xi (q))$ and $\hat{\Xi }(p)=(\hat{\pi }^2_0(p),\hat{\xi
}(p))$, then
$$\hat{\xi }(p)=\xi (\hat{\pi }^2_2(p))\sigma (p).$$

Using Lemma \ref{rbsl5} we have that $\sigma $, and $\Sigma $, are
differentiable.

Let $p\in (\hat{\pi }^2_2)^{-1}(U)$. Since $(\hat{\pi
}^2_2)^{-1}(\hat{\pi }^2_2(p))=\hat{\pi }^2_2(p)A_2(n)$ (see Lemma
\ref{rbsl2}), then we can put $p=\hat{\pi }^2_2(p)(I,h)$, with
$h\in A_2(n)$. Therefore
$$\xi (\hat{\pi }^2_2(p))(I,h)=\hat{\xi }(\hat{\pi
}^2_2(p))(I,h)=\hat{\xi }(\hat{\pi }^2_2(p)(I,h))=\hat{\xi }(p).$$

Comparing the definition of $\sigma $ with the former expression
and keeping in mind Lemma \ref{grol1}, iii), we obtain
$$\sigma (p)=(I,h).$$

Moreover, $(\hat{\pi }^2_2)^{-1}(q)=qA_2(n)$ for every $q\in
{\mathcal F}^2M$; therefore $\Sigma $ is a bijection and, as a
consequence, $\Sigma $ is a diffeomorphism.

Finally we obtain that the map $\sigma $ verifies $\sigma
(p(I,\tilde{h}))=\sigma (p)(I,\tilde{h})$ for every $\tilde{h}\in
A_2(n)$ and $p=\hat{\pi }^2_2(p)(I,h)\in (\hat{\pi
}^2_2)^{-1}(U)$. Of course, we should see that $\sigma
(p(I,\tilde{h}))=(I,h)(I,\tilde{h})$, but this is immediate from
the very definition of $\sigma $ and from the fact that
$$\xi (\hat{\pi }^2_2(p(I,\tilde{h})))(I,h)(I,\tilde{h})=\xi (\hat{\pi
}^2_2(p))(I,h)(I,\tilde{h})=\hat{\xi }(p)(I,\tilde{h})=\hat{\xi
}(p(I,\tilde{h})).$$

This completes the proof.\\

In~\cite{paula} it is proved that $\tilde{\mathcal F}^2M$ can be
endowed with a $G^1(n)$-principal bundle structure over
$\hat{\mathcal F}^2M$. Following the lines in the proof of Theorem
5.2 in~\cite{paula}, straightforward calculations show that the
projection $\pi \colon \tilde{\mathcal F}^2M \rightarrow
\hat{\mathcal F}^2M$ defined there is given in local coordinates
by $\pi
(x^{i},x^{i}_j,y^{i}_j,x^{kl}_j)=(x^{i},x^{i}_j,x^{kl}_rx^r_j)$ (
in other words, $\pi (x,a,b,f)=(x,a,f (I, a))$), and that the
action on the right of the group $G^1(n)$ on $\tilde{\mathcal
F}^2M$ is expressed, identifying $l\in G^1(n)$ with $(I,l,0)\in
\tilde{G}^2(n)$, as $(x,a,b,f)(I,l,0)=(x,a,bl,f (I, l))$. This
action is, in fact, the restriction of the action on the right of
$\tilde{G}^2(n)$ on $\tilde{\mathcal F}^2M$.

Now we must identify $h\in A_2(n)$ with $(I,I,h)\in
\tilde{G}^2(n)$. Let us consider the subgroup of $\tilde{G}^2(n)$
given by the semi-direct product $\tilde{G}^2_2(n)=G^1(n)\times
A_2(n)$, where the product law is given by
$$(I,l,h)(I,l^{\prime },h^{\prime })=(I,ll^{\prime
},h^{\prime}+h(I,l^{\prime})).$$

Notice to warnings: first, if $(I,l,0)\in G^1(n)$ and $(I,I,h)\in
A_2(n)$, then we obtain $(I,l,0)(I,I,h)=(I,l,h)$; second, this
group $\tilde{G}^2_2(n)$ has nothing to do with the group
$G^1(n)\times A_2(n)$ of Lemma \ref{grol1}, ii), which is a
subgroup of $\hat{G}^2(n)$.

\begin{teo}
\label{rbst2} $\tilde{\mathcal F}^2M$ is a principal bundle over
${\mathcal F}^2M$ with structure group $\tilde{G}^2_2(n)$ and
projection $\tilde{\pi }^2_2=\hat{\pi }^2_2 \circ \pi$.
\end{teo}

{\bf Proof:}

$ $

The action on the right of the group $\tilde{G}^2_2(n)$ on
$\tilde{\mathcal F}^2M$ is the restriction of the action on the
right of $\tilde{G}^2(n)$ on $\tilde{\mathcal F}^2M$; so, such
action is free. This action can be described as
$$(x,a,b,f)(I,l,h)=(x,a,bl,a\circ h + f(I,l)),$$

\noindent an expression that can be obtained as well as
$((x,a,b,f)(I,l,0))(I,I,h)$. Let us also point out that the group
law in $\tilde{G}^2_2(n)$ is identical to the law of the group
$\tilde{G}^2_1(n)$ introduced at the end of Section \ref{pre}, so
the fact that $p ((I,l,h)(I,l^{\prime},h^{\prime }))=(p
(I,l,h))(I,l^{\prime},h^{\prime })$ for every $p\in
\tilde{\mathcal F}^2M$ and $(I,l,h)$, $(I,l^{\prime},h^{\prime
})\in \tilde{G}^2_2(n)$ follows from the fact that
$\tilde{\mathcal F}^2M$ is a $\tilde{G}^2_1(n)$-principal bundle
over ${\mathcal F}M$.

The rest of the proof is immediate.\\

It is also obvious, from the definition of $\tilde{\pi }^2_2$
itself, that the projection $\tilde{\pi }^2_2$ satisfy the
appropriate relations of commutativity with all the other
projections.

\section{Conclusion}
\label{con}

\noindent Theorems \ref{rbst1} and \ref{rbst2} provide the bundle
structures missing in~\cite{paula} and conclude the relations of
existence of principal bundle structures among the five manifolds
$\tilde{\mathcal F}^2M$, $\hat{\mathcal F}^2M$, ${\mathcal F}^2M$,
${\mathcal F}M$ and $M$ we are dealing with. We collect all of
them in the following picture, where we have pointed out at the
right side of the arrows the corresponding structure group:

$$\begin{array}{cccccccccc}
 &  &  &  \tilde{\mathcal F}^2M &  & \hat{\mathcal F}^2M &   &  {\mathcal
F}^2M &  &  {\mathcal F}M \\
 \stackrel{descending \; to \;}{ \scriptstyle the \; first \; level} & & &
\makebox[0pt][l]{\;\;\;\;$\scriptstyle \tilde{G}^2(n)$} \downarrow
& & \downarrow
 \makebox[0pt][l]{$\scriptstyle \hat{G}^2(n)$} & & \downarrow
 \makebox[0pt][l]{$\scriptstyle G^2(n)$} & & \downarrow
 \makebox[0pt][l]{$\scriptstyle G^1(n)$} \\
 &  &  & M & & M &   &  M &  &  M
 \end{array}$$

 \vspace{0.5cm}

 $$\begin{array}{cccccccccccc}
 &  & &  &  &  \tilde{\mathcal F}^2M &  & \hat{\mathcal F}^2M &   &  {\mathcal
F}^2M  &  & \\
 \stackrel{descending \; to \;}{ \scriptstyle the \; second \; level} &  & &  &  &
\makebox[0pt][l]{\;\;\;\;$\scriptstyle \tilde{G}^2_1(n)$}
\downarrow & & \downarrow
 \makebox[0pt][l]{$\scriptstyle L_2(n)$} & & \downarrow
 \makebox[0pt][l]{$\scriptstyle S_2(n)$}  &  & \\
 &  & &  &  & {\mathcal F}M & & {\mathcal F}M &   &  {\mathcal F}M
 &  &
 \end{array}$$

 \vspace{0.5cm}

 $$\begin{array}{ccccccccc}
 &  & &  &  &  \tilde{\mathcal F}^2M &  & \hat{\mathcal F}^2M &    \\
 \stackrel{descending \; to \;}{ \scriptstyle the \; third \; level}  & \;\;\;\;\;\;\;\;\;& &  &  &
\makebox[0pt][l]{\;\;\;$\scriptstyle \tilde{G}^2_2(n)$} \downarrow
& & \downarrow
 \makebox[0pt][l]{$\scriptstyle A_2(n)$} & \;\;\;\;\;\;\;\; \\
 &  & &  &  & {\mathcal F}^2M & & {\mathcal F}^2M &
 \end{array}$$

 \vspace{0.5cm}

 $$\begin{array}{cccccccccccccc}
 &  &  &  &  &  &  \tilde{\mathcal F}^2M &  &  & \\
 \stackrel{descending \; to \;}{ \scriptstyle the \; fourth \; level}  & \;\;\;\;\; &  & \;\;\;\;\;\; &  &  &
\makebox[0pt][l]{\;\;\;\;$\scriptstyle G^1(n)$} \downarrow & \;\;\;\;\;\; & \;\;\;\;\;\; & \\
 &  &  &  &  &  &  \hat{\mathcal F}^2M &  &  &
 \end{array}$$

\vspace{0.5cm}

For a detailed description of the principal bundle structures not
treated here, see~\cite{paula} and references therein.

{\bf Acknowledgements}

This work has been partially supported by MINECO MTM2013-42870-P,
the European project IRSES-project "Geomech-246981" and the ICMAT
Severo Ochoa project SEV-2011-0087.

\end{document}